\documentclass[final]{siamltex}

\usepackage{amsfonts}
\usepackage{amssymb,amsmath}
\usepackage[mathscr]{eucal}
\usepackage{graphicx}
\usepackage{times}
\usepackage{float}

\title{Alternative Orthogonal Rational Functions\\ on a Half Line}

\author{Vladimir S. Chelyshkov \thanks{ University of Pikeville, 147 Sycamore Street. Pikeville, Kentucky 41501 \newline
\hspace{10pt}({\tt VolodymyrChelyshkov@upike.edu}).
}}

\begin{document}

\maketitle

\begin{abstract}
System of alternatively orthogonalized rational functions of Jacobi type on the half line $[1, \infty)$ is defined and its properties are established.
Three subsystems of proper and mixed systems of rational functions with nice properties are presented. 
\end{abstract}



\pagestyle{myheadings}
\thispagestyle{plain}
\markboth{V.~S.~CHELYSHKOV}{ORTHOGONAL RATIONAL FUNCTIONS}

\section{Introduction}

Interest in approximation on a semi-infinite interval  by algebraically mapped classical orthogonal polynomials on a finite interval originates from \cite{Grosch}, \cite{Boyd}, and such an approach was later developed by many authors (see, for example \cite{Ben2}, \cite{Ben1}).

By following this path, we present system of alternative orthogonal rational functions of Jacobi type defined on a half line. We find that the functions possess all the same properties as the alternative orthogonal polynomials do. Also, we show that algorithmic capabilities of approximation by  constructed (both proper and mixed) subsystems of rational functions are equivalent to those of the classical orthogonal polynomials.

\section{Bidirectional orthogonalization and the Jacobi polynomials}

Procedure of bidirectional orthogonalization of a sequence of functions and alternative orthogonal polynomials were introduced and studied in
\cite{VC1}, \cite{VC4}, \cite{VC6}.

Following \cite{VC6}, let us consider two systems polynomials
\begin{equation}
\{{\cal P}_{nk}^{(\alpha,\beta)}(x)\}_{k=n}^{0},\,\, \alpha>-2,\,\, \beta>-1, \quad
 \mbox{and} \quad \{P_{nk}^{(\alpha,\beta)}(x)\}_{k=n}^{\infty},\,\, \alpha,\beta>-1
\label{2sys}
\end{equation}
that are orthogonal with respect to the weight function $w(x;\alpha,\beta)=x^{\alpha}(1-~x)^{\beta}$ on the interval $[0,1]$.  The systems are constructed by applying inverse and direct  orthogonalization algorithms beginning with $k=n$. The polynomials in (\ref{2sys}) are related by identities 
\begin{equation}
{\cal P}_{nk}^{(\alpha,\beta)}(x)=P_{-n,-k}^{(\gamma,\beta)}(x^{-1})=x^{-1}P_{-(n+1),-(k+1)}^{(\gamma+2,\beta)}(x^{-1})= ...,
\label{j6}
\end{equation}
\[
\mbox{with} \quad 0\leq k \leq n \quad \mbox{and} \quad \gamma=-\alpha-\beta-2.
\]
Also,
\begin{equation}
{\cal P}_{nk}^{(\alpha,\beta)}(x)=x^{k}{P}_{0,n-k}^{(\alpha+2k+1,\beta)}(x), \quad  0\leq k \leq n,
\label{jj3}
\end{equation}
and
\begin{equation}
\label{j3}
P_{nk}^{(\alpha,\beta)} (x)=x^nP_{0,k-n}^{(\alpha+2n,\beta)} (x) \quad k\ge n,
\end{equation}
where $P_{0m}^{(\alpha,\beta)} (x)$ are the Jacobi polynomials on the interval $[0,1]$. 

Identities (\ref{j6}) - (\ref{j3}) were used  in \cite{VC6} to derive classical properties  of the alternative Jacobi polynomials ${\cal P}_{nk}^{(\alpha,\beta)}(x)$.

\section{Definitions and properties} 

We obtain alternative orthogonal rational functions on the interval $[1,\infty)$ making use of key identities (\ref{j6}) and other properties of the polynomials ${\cal P}_{nk}^{(\alpha,\beta)}(x)$ and $P_{nk}^{(\alpha,\beta)}(x)$ established in \cite{VC6}.

\paragraph{ First definition and orthogonality} 

We introduce 
\[
w(x;\alpha,\beta)=x^\alpha |1-x|^\beta
\]
as a unification formula for the weight functions on the intervals $[0,1]$ and $[1,\infty)$ and define alternative orthogonal rational functions as follows
\[
{\cal R}_{nk}^{(\gamma,\beta)}(x)={\cal P}_{nk}^{(\gamma,\beta)}(x^{-1}),
\label{r0}
\]
\[
0\le k\le n, \quad \beta>-1\quad \mbox{and}\quad \gamma+\beta<0. 
\] 
The functions   ${\cal R}_{nk}^{(\gamma,\beta)}(x)$ satisfy the orthogonality relation
\begin{equation} 
\int_1^\infty (x-1)^\beta x^\gamma  {\cal R}_{nk}^{(\gamma,\beta)}(x){\cal R}_{nl}^{(\gamma,\beta)}(x)
{\rm d}x=r_{nk}^{(\gamma,\beta)}\delta_{kl},
\label{r3}
\end{equation} 
where 
\[
r_{nk}^{(\gamma,\beta)}=\frac{1}{2k-\gamma-\beta-1}\frac{\Gamma(n+k-\gamma-\beta)}{\Gamma(n+k-\gamma)}\frac{\Gamma(n-k+\beta+1)}{(n-k)!}.
\]

\paragraph{Invariance}
For $p \in {\mathbb N}$
 orthogonality  (\ref{r3})  holds the invariance
\[
{r}_{n+p,k+p}^{(\gamma+2p,\beta)}={r}_{nk}^{(\gamma,\beta)},
\]
and 
\[
{\cal R}_{nk}^{(\gamma,\beta)}(x)={\cal R}_{n+p,k+p}^{(\gamma+2p,\beta)}(x)=x^{-p}{\cal R}_{nk}^{(\gamma+2p,\beta)}(x).
\]

\paragraph{Integral}
In addition, one can evaluate
\[
\hspace{-66pt} \int_1^{\infty}(1-x)^\beta x^\gamma {\cal R}_{nk}^{(\gamma,\beta)}(x) {\rm d} x
\]
\[
=\frac{\Gamma(k-\gamma-\beta-1)}{k!}\frac{\Gamma(\beta+n-k+1)}{(n-k)!}\frac{n!}{\Gamma(n-\gamma)}.
\]

\paragraph{Recurrence Relations}
With
\[
{\cal R}_{nn}^{(\gamma,\beta)}(x)=x^{-n},
\]
\[
{\cal R}_{n,n-1}^{(\gamma,\beta)}(x)=-(\gamma+ \beta -2n+2)x^{-(n-1)}+(\gamma - 2n+1)x^{-n}
\]
the functions ${\cal R}_{nk}^{(\gamma,\beta)}(x)$ satisfy the three-term recurrence relation 

\[
(n-k+1)(\gamma + \beta-n-k+1)(\gamma + \beta+2k){\cal R}_{n,k-1}^{(\gamma,\beta)}(x)
\]
\[
=(\gamma + \beta-2k+1)[-(\gamma+ \beta - 2k + 2)(\gamma + \beta-2k)x
\]
\[ 
- (\gamma+ \beta - 2n)(\gamma - 2k+1)-2(n-k)(n-k+1)]{\cal R}_{nk}^{(\gamma,\beta)}(x)
\]
\[
-(\gamma - n - k)(\beta + n - k)(\gamma+ \beta - 2k + 2){\cal R}_{n,k+1}^{(\gamma,\beta)}(x).
\]

\paragraph{Second definition}
Identities (\ref{j6}) can also be represented as ( cf.  \cite{VC1})
\[
{\cal P}_{nk}^{(\gamma,\beta)}(x^{-1})= P_{-n,-k}^{(\gamma,\beta)}(x)=xP_{-(n+1),-(k+1)}^{(\gamma+2,\beta)}(x),\quad k=n,n-1 ,..., 0,
\]
which results in definition
\begin{equation}
{\cal R}_{nk}^{(\gamma,\beta)}(x)= P_{-n,-k}^{(\gamma,\beta)}(x).
\label{r2}
\end{equation}
We employ (\ref{r2}) for deriving results given below.

\paragraph{ Rodrigues' type formula and integral representation for the rational functions}
\[
{\cal R}_{nk}^{(\gamma,\beta)}(x)=\frac{(-1)^{n-k}x^{n-\gamma}(x-1)^{-\beta}}{(n-k)!}\frac{{\mbox{d}}^{n-k}}{{\mbox{d}}x^{n-k}}(x^{\gamma-n-k}(x-1)^{\beta+n-k})
\]
with $k=n, n-1, ..., 0$, and 

\[
{\cal R}_{nk}^{(\gamma,\beta)}(x)=\frac{(-1)^{n-k}}{2\pi i}\frac{x^{n-\gamma}}{(x-1)^{\beta}}\int\limits_{C}
{\frac{z^{\gamma-n-k}(z-1)^{\beta+n-k}}{(z-x)^{n-k+1}}}{\rm d}z.
\]
Here $C$ is a closed contour encircling the point $z=x$.

\paragraph{Differential-difference relations} 
\[
(2k+\gamma+2)x(x-1)\frac{\rm d}{{\rm d}x}{\cal R}_{nk}^{(\gamma,\beta)}(x)=[n(n+\gamma+\beta+2) +k(k-\beta )
\]
\[
- k(2k+\gamma+2)x]{\cal R}_{nk}^{(\gamma,\beta)}(x)+(n+k+\gamma+\beta+2)(n-k+\beta){\cal R}_{n,k+1}^{(\gamma,\beta)}(x)
\]
and
\[
(2k+\gamma)x(x-1)\frac{\mbox{d}}{\mbox{d}x}{\cal R}_{nk}^{(\gamma,\beta)}(x)=[-n(n+\gamma+\beta+2) - (k+\gamma+\beta-1)(k+\gamma+\beta+1)
\]
\[
+(2k+\gamma)(k+\gamma+1)x]{\cal R}_{nk}^{(\gamma,\beta)}(x)- (n-k+1)(n+k+\gamma-1){\cal R}_{n,k-1}^{(\gamma,\beta)}(x).
\]

\paragraph{Differentiation formula} 
\[
\frac{\rm d}{{\rm d}x}{\cal R}_{nk}^{(\gamma,\beta)}(x)+\frac{n}{x}{\cal R}_{nk}^{(\gamma,\beta)}(x)=(n + k+\gamma+1){\cal R}_{n,k+1}^{(\gamma+1,\beta+1)}(x).
\]

\begin{figure}[H]
\centerline{\includegraphics[height=36mm, width=60mm]{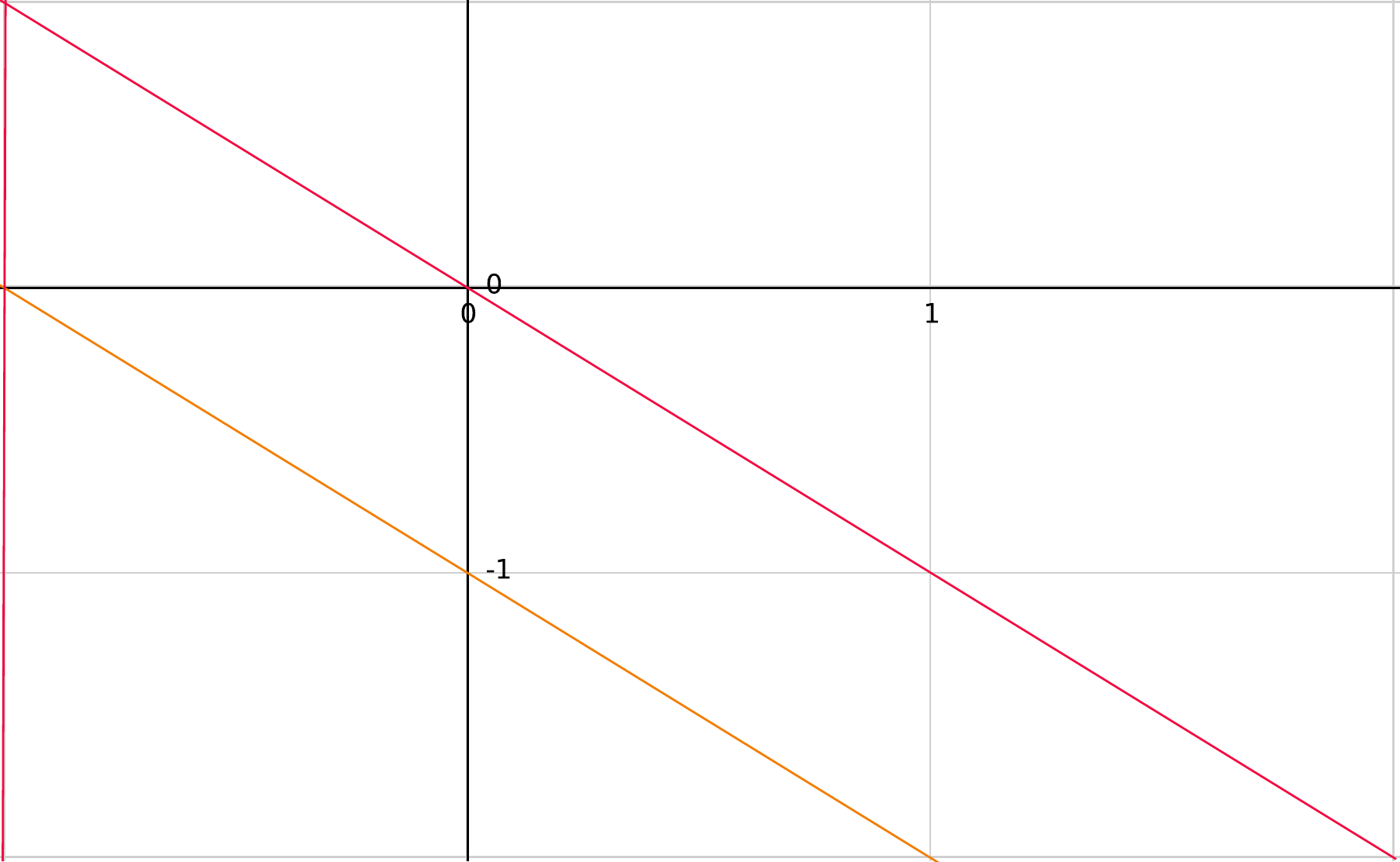}}
\caption{$(\beta,\gamma)$-parametric domain with the  oblique  marginal strip.}
\normalsize
\label{domain}
\end{figure}

\paragraph{Differential equation}
The function $y(x)={\cal R}_{nk}^{(\gamma,\beta)}(x)$ is a solution to the differential equation
\[
x^2(x-1)y^{\prime\prime}+x((\gamma+\beta+2) x - \gamma-1)y^{\prime}
-[k(k-\gamma-\beta-1)x-n(n-\gamma)]y=0
\]
with $(\gamma, \beta)$ in the domain specified below (see Fig. \ref{domain}).

\section{Gaussian quadratures and discrete orthogonality}
Non-trivial zeros of  ${\cal R}_{n,k-1}^{(\gamma,\beta)}(x)$  provide abscissas of alternative $(\gamma,\beta)$-Gaussian quadratures that are exact for $x^{-j}$ with 
\begin{equation}
2k-1 \le j \le 2n, \quad 0<k \le n.
\label{g1}
\end{equation}
Inequalities (\ref{g1}) grant discrete orthogonality, but a chosen quadrature is exact for $k \le j \le 2n$,  iff $k=1$. We call the quadrature with $k=1$  $the$ alternative $(\gamma,\beta)$-Gaussian quadrature for rational functions. 

Similarly, one may also find that orthogonal rational functions defined by direct algorithm of orthogonalization as $R^{(\gamma,\beta)}_{nk}(x)=P_{nk}^{(\gamma,\beta)}(x^{-1})$ generate Gauss-type quadrature with the zeros of $R^{(\gamma+2n,\beta)}_{0,k-n}(x)$. The quadrature is exact for $x^{-j}$ with 
\[
2n \le j \le 2k-1, \quad k>n.
\]
It becomes the $(\gamma,\beta)$-Gaussian quadrature for rational functions iff $n=0$.

\section{Proper alternative orthogonal rational functions}
Let $n \in  {\mathbb N}$ and
\[
\mbox{\boldmath ${\cal R}$}_{n}^{(\gamma,\beta)}(x)=\{{\cal R}_{nk}^{(\gamma,\beta)}(x)\}_{k=n}^1 
\]
be a system of proper rational functions.  
This orthogonal system is of special interest inasmuch as its terms are quadratically integrable functions:
\begin{equation}
\int_1^{\infty} \left({\cal R}_{nk}^{(\gamma,\beta)}(x)\right) ^{2} {\rm d}x <\infty.
\label{proper}
\end{equation}
The associated orthogonal function ${\cal R}_{n0}^{(\gamma,\beta)}(x)$ does not satisfy (\ref{proper}),  but the zeros of  ${\cal R}_{n0}^{(\alpha,\beta)}(x)$ are abscissas of the alternative $(\alpha,\beta)$-Gaussian quadrature for proper rational functions. Abscissas and weights of the quadrature can be  recalculated from the standard $(\alpha,\beta)$-Gauss quadrature over a closed interval. We also find that $
\mbox{\boldmath ${\cal R}$}_{n}^{(\gamma,\beta)}(x)$  possesses discrete orthogonality property with the abscissas  chosen.

All of this makes the set 
\[
\left(
\mbox{\boldmath ${\cal R}$}_{n}^{(\gamma,\beta)}\cup{\cal R}_{n0}^{(\gamma,\beta)}\right)(x)
\]
a system of its own kind.

\section{Singular orthogonality and marginal systems} 
Let
\vspace{-10pt}
\begin{equation} 
\gamma^{\prime}\in\{\gamma\,\,|\,\,\,-1\leq\gamma+\beta <0,\,\,\,\, -1<\beta\}. 
\label{sng}
\end{equation}
In this exceptional case the improper rational function ${\cal R}_{n0}^{(\gamma^{\prime},\beta)}(x)$ is not square-integrable with respect to the weight function $w(x;\gamma^{\prime},\beta)$ and therefore non-normalizable. However,  it is consistent with orthogonality relations (\ref{r3})  for $l>~0$. Thus,  ${\cal R}_{n0}^{(\gamma^{\prime},\beta)}(x)$   is a singular term of the set
\[
 \left(
\mbox{\boldmath ${\cal R}$}_{n}^{(\gamma^{\prime},\beta)}\cup{\cal R}_{n0}^{(\gamma^{\prime},\beta)}\right)(x)
\]
which we call a marginal system. Below we describe two noteworthy orthogonal systems of such a kind.

\section{A-kind rational functions}
Considering (\ref{sng}) we find that functions ${\cal R}_{n0}^{\mathscr A}(x) ={\cal R}_{n0}^{(-1,0)}(x)$
are the shifted Legendre polynomials of the reciprocal $1/x$, and we  count ${\cal R}_{n0}^{\mathscr A}(x)$ as the associated term for the system of proper rational functions
\[
\mbox{\boldmath ${\cal R}$}_{n}^{\mathscr A}(x)=\{{\cal R}_{nk}^{\mathscr A}(x)\}_{k=n}^1, \quad {\cal R}_{nk}^{\mathscr A}(x) = {\cal R}_{nk}^{(-1,0)}(x).
\]
The functions ${\cal R}_{nk}^{\mathscr A}(x)$ obey the orthogonality relations
\[
\int_1^{\infty}\frac{1}{x}
{\cal R}_{nk}^{\mathscr A}(x)
{\cal R}_{nk}^{\mathscr A}(x){\rm d}x=\frac{\delta_{kl}}{k+l},\quad k=0,1,...,n,\quad l=1,2,...,n
\]
\begin{figure}[H]
\centerline{\includegraphics[height=36mm, width=60mm]{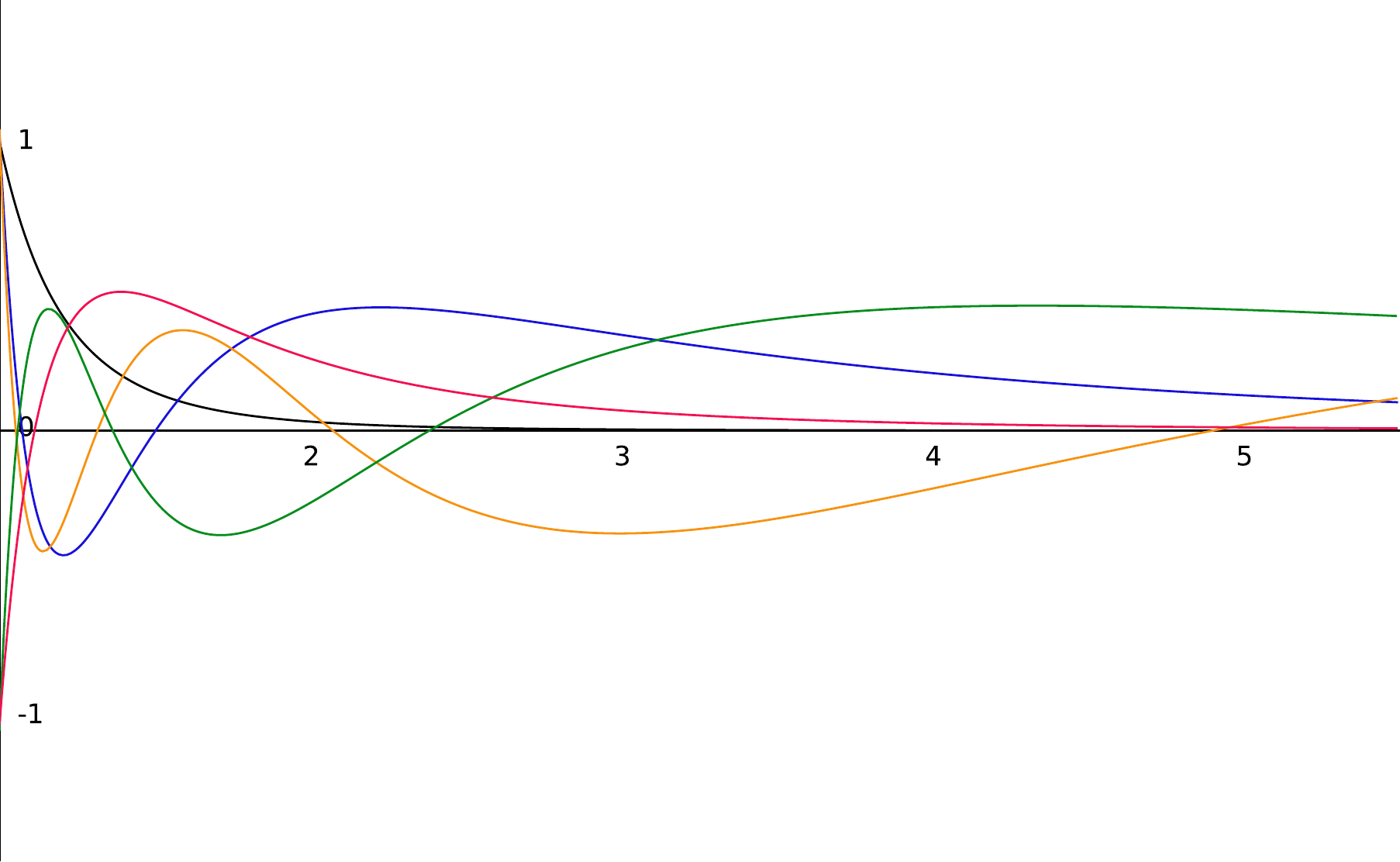}}
\caption{${\cal R}_{nk}^{\mathscr A}(x)$: $n=5, k=1-5$.}
\normalsize
\end{figure}
\hspace{-20pt}
and 
\[
\int_1^{\infty}\frac{1}{x}
{\cal R}_{nk}^{\mathscr A}(x){\rm d}x=\frac{1}{k},\quad k=1,...,n.
\]
They can be calculated using the expansion
\[
{\cal R}_{nk}^{\mathscr A}(x)=\sum\limits_{j=0}^{n-k} {(-1)^j} \left( {\begin{array}{c}
n-k\\
\mbox{ }j \\
\end{array}} \right)\left({\begin{array}{c}
n+k+j\\
\mbox{ }n-k\\
\end{array}}\right)x^{-(k+j)},\mbox{ }k=0,1,...,n,
\]
or the three-term recurrence relation
\[
{\cal R}_{nn}^{\mathscr A}(x)=x^{-n},\quad {\cal R}_{n,n-1}^{\mathscr A}(x)=(2n-1)x^{-(n-1)}-2nx^{-n},
\]
\[
(2k+1)(n+k)(n-k+1){\cal R}_{n,k-1}^{\mathscr A}(x)
\]
\[
=2k[(2k-1)(2k+1)x-2(n^{2}+k^{2}+n)]{\cal R}_{nk}^{\mathscr A}(x)
\]
\[
-(2k-1)(n-k)(n+k+1){\cal R}_{n,k+1}^{\mathscr A}(x).
\vspace{8pt}
\]
We also find that $y(x)={\cal R}_{nk}^{\mathscr A}(x)$    is a solution to the differential equation

\[
x^2(x-1)y^{\prime\prime}+x^2 y^{\prime}-[k^2 x-n(n+1)]y=0.
\]
The set
\[
\left(
\mbox{\boldmath ${\cal R}$}_{n}^{\mathscr A}\cup{\cal R}_{n0}^{\mathscr A}\right)(x)
\]
provides regular Gauss-Rational quadrature over a half line, which is exact for for $x^{-j}, \,\, 1 \le j \le 2n$ with the weight $1/x$,  and discrete orthogonality for $\mbox{\boldmath ${\cal R}$}_{n}^{\mathscr A}(x)$.

\section{T-kind rational functions}
\label{sec: 6}
Let notation $(a)_n$ represents falling factorial. Considering (\ref{sng}) we find that 
\[ 
{\cal R}_{n0}^{\mathscr T}(x) = n!/(n-\tfrac{1}{2})_{n}{\cal R}_{n0}^{(0,-1/2)}(x)
\]
are the shifted Chebyshev polynomials of the reciprocal $1/x$, and we  count ${\cal R}_{n0}^{\mathscr T}(x)$ as the associated term for the system of proper rational functions
\[
\mbox{\boldmath ${\cal R}$}_{n}^{\mathscr T}(x)=\{{\cal R}_{nk}^{\mathscr T}(x)\}_{k=n}^1
\]
\[
{\cal R}_{nk}^{\mathscr T}(x) = c_{nk}{\cal R}_{nk}^{(0,-1/2)}(x),\quad c_{nk}=(n-k)!/(n-k-1/2)_{n-k},
\]
$n \in  {\mathbb N}$. 
\begin{figure}[htbp]
\centerline{\includegraphics[height=36mm, width=60mm]{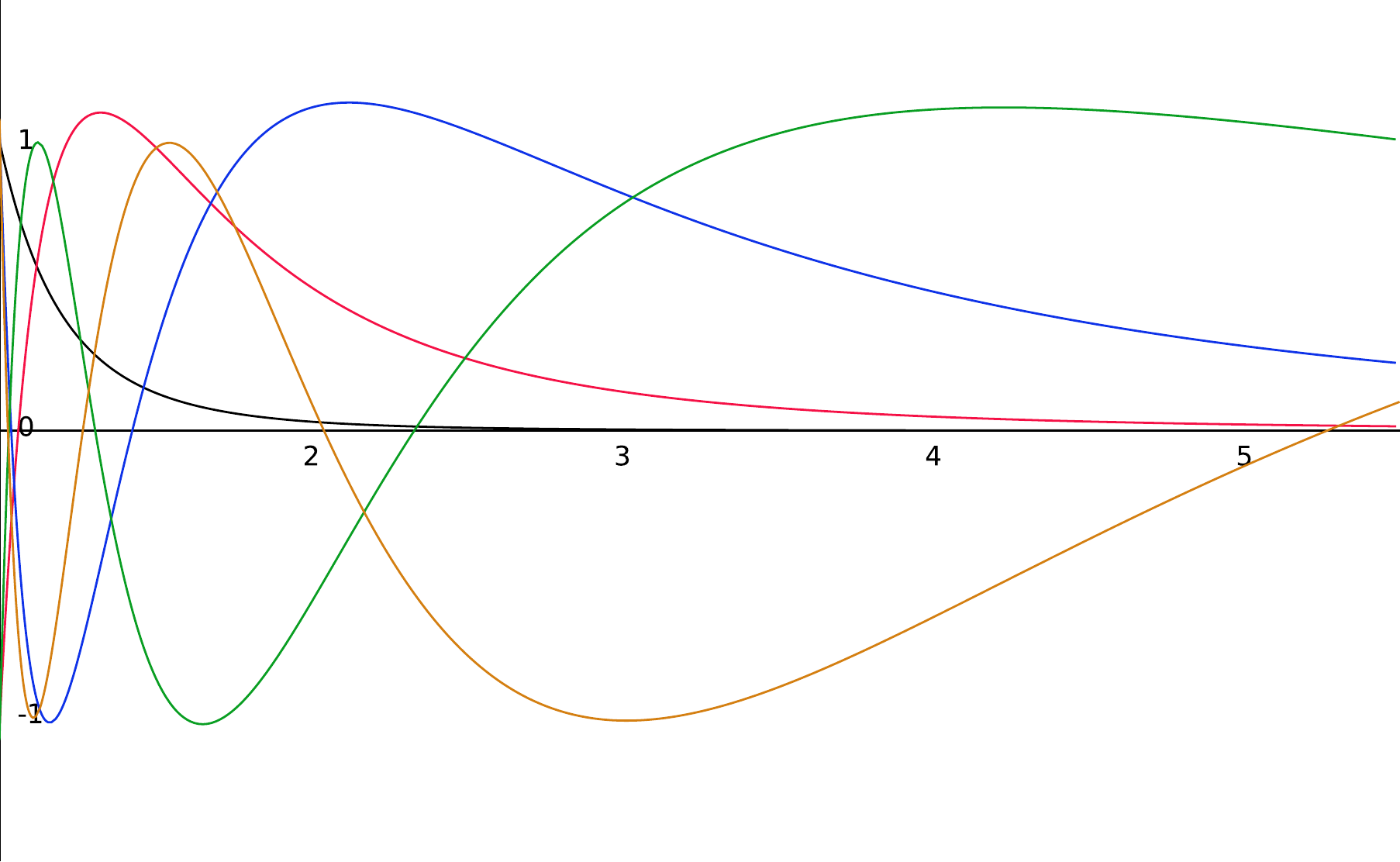}}
\caption{${\cal R}_{nk}^{\mathscr T}(x)$: $n=5, k=1-5$.}
\normalsize
\end{figure}
The functions ${\cal R}_{nk}^{\mathscr T}(x)$  obey the orthogonality relation
\[
\int_1^{\infty}\frac{{\cal R}_{nk}^{\mathscr T}(x){\cal R}_{nk}^{\mathscr T}(x)
}{\sqrt{x-1}}{\rm d}x=\frac{1}{2(k+l)-1}\frac{(2n+k+l-1)!!}{(2n+k+l-2)!!}\frac{(2n-k-l)!!}{(2n-k-l-1)!!}\pi\delta_{kl},
\]
\[
k=0,1,...,n,\quad l=1,2,...,n,
\]
 and
\[
\int_1^{\infty}\frac{{\cal R}_{nk}^{\mathscr T}(x)}   
{\sqrt{x-1}}{\rm d}x=\frac{(2k-3)!!}{(2k)!!}2n\pi,\quad k=1,2...n.
\]
They can be calculated using the three-term recurrence relation
\[
{\cal R}_{nn}^{\mathscr T}(x)=x^{-n},\quad {\cal R}_{n,n-1}^{\mathscr T}(x)=(4n-3)x^{-(n-1)}-(4n-2)x^{-n},
\]
\[
[2(n+k)-1][2(n-k)+1](4k+1){\cal R}_{n,k-1}^{\mathscr T}(x)
\]
\[
=(4k-1)[(4k-3)(4k+1)x-2(4n^{2}+4k^{2}-2k-1)]{\cal R}_{nk}^{\mathscr T}(x)
\]
\[
-4(n-k)(n+k)(4k-3){\cal R}_{n,k+1}^{\mathscr T}(x).
\]
We also find that $y(x)={\cal R}_{nk}^{\mathscr T}(x)$ is a solution to the differential equation

\[
x^2(x-1)y^{\prime\prime}+x(\tfrac{3}{2} x  - 1) y^{\prime}-[k(k-\tfrac{1}{2}) x - n^2]y=0.
\]
The set
\[
\left(
\mbox{\boldmath ${\cal R}$}_{n}^{\mathscr T}\cup{\cal R}_{n0}^{\mathscr T}\right)(x)
\]
provides Chebyshev-Gauss-Rational quadrature over a half line and discrete orthogonality for $\mbox{\boldmath ${\cal R}$}_{n}^{\mathscr T}(x)$. 

In addition, a fast Fourier transform can be employed for  computations.

\section{Mixed systems of alternative orthogonal rational functions}

Let 
\[
\gamma+\beta<-1,\quad \beta>-1
\] 
and 
\[
\mbox{\boldmath $\widetilde{{\cal R}}$}_{n}^{(\gamma,\beta)}(x)=\{{\cal R}_{nk}^{(\gamma,\beta)}(x)\}_{k=n}^0
\]
be a mixed system, i.e., the system that contains an  improper rational function with $k=0$. Such a system may be of interest for approximation of functions that have a finite value at infinity. 

We introduce alternative Legendre rational functions
\[
 \mbox{\boldmath $\widetilde{{\cal R}}$}_{n}(x)=\mbox{\boldmath $\widetilde{{\cal R}}$}_{n}^{(-2,0)}(x),
\] 
considering them as the simplest subsystem of $\mbox{\boldmath $\widetilde{{\cal R}}$}_{n}^{(\gamma,\beta)}(x)$  with nice properties. Functions
\[
 \mbox{\boldmath $\widetilde{{\cal R}}$}_{n}(x)=\{{\cal R}_{nk}(x)\}_{k=n}^0
\] 
obey the orthogonality relation
\[
\int_1^{\infty}\frac{1}{x^{2}}
{\cal R}_{nk}(x)
{\cal R}_{nk}(x){\rm d}x=\frac{\delta_{kl}}{k+l+1},\quad k,l=0,1,...,n,
\]
and their properties directly follow from the properties of alternative Legendre polynomials described in \cite{VC4}. The system
$\mbox{\boldmath $\widetilde{{\cal R}}$}_{n}(x)$
provides Radau-Rational quadrature over a half line with one of the abscissas at infinity; also, it possesses discrete orthogonality for $\{{\cal R}_{nk}(x)\}_{k=n}^1$. Abscissas and weights of the quadrature can be recalculated from the standard Radau quadrature over a closed interval.


\section{Conclusion} 

Orthogonal functions presented in this preprint preserve hierarchy of algebraic decay of the sequence $x^{-j}$, $j=n, n-1, ..., 1$ at infinity. This property 
may be considered as an advantage for approximation on a half line, in particular if the asymptotic of a function is not known exactly. Weak formulation for problems involving differentiation may be required. Also, the systems may be of interest for studying energy cascades in nonlinear problems.


\begin{thebibliography}{}

\bibitem{Boyd}
{\sc J.P. Boyd},  {\em Spectral methods using rational basis functions on an infinite interval}, J. Comput. Phys., Vol. 69, No 1, (1987), pp.~112-142. 

\bibitem{VC1}
{\sc  V.~S. Chelyshkov}, {\em  Sequences of exponential polynomials, which are orthogonal on the semi-axis}, 
Dokl. Akad. Nauk Ukr. SSR,  Ser. A, No 1, (1987) pp.~14 --17.  

\bibitem{VC4}
{\sc  V.~S. Chelyshkov}, {\em Alternative orthogonal polynomials and quadratures}, Electron. Trans. Numer. Anal., Vol. 25, (2006), pp.~17-26. 

\bibitem{VC6} 
{\sc V.~S. Chelyshkov}, {\em Alternative Jacobi polynomials and orthogonal exponentials}, arXiv:1105.1838.

\bibitem{Grosch}
{\sc C.E. Grosch, S.A.  Orszag}, {\em Numerical solution of problems in unbounded regions: coordinate transforms}, J. Comput. Phys., Vol. 25, No 3, (1977), pp. 273-295. 

\bibitem{Ben2}
{\sc Yong-gang Yi, Ben-yu Guo}, {\em  Generalized Jacobi rational spectral method on the half line},  Adv. Comput. Math., Vol. 37, No 1, (2012), pp. 1-37.

\bibitem{Ben1}
{\sc Zhong-Qing Wang, Ben-Yu Guo},  {\em Jacobi rational approximation and spectral methods for differential equations of degenerate type}, Math. Comp., Vol. 77, No 262, (2008), pp. 883-907.

\end{thebibliography}
\end{document}